\documentclass[11pt]{article}
\usepackage{amsmath,amssymb,latexsym,epsfig,subfigure,multirow,graphicx,longtable,rotating}

\newcommand{\q}{\quad}

\newcommand{\na}{\nabla}
\newcommand{\hs}{\hspace}

\newcommand{\txf}{\textbf}

\newcommand{\de}{\delta}

\newcommand{\veps}{\varepsilon}

\newcommand{\ka}{\kappa}

\newcommand{\ph}{\phi}

\newcommand{\Ga}{\Gamma}

\newcommand{\Om}{\Omega}

\newcommand{\mb}{\mathbb}

\newcommand{\iy}{\infty}
\newcommand{\fr}{\frac}
\newcommand{\pa}{\partial}

\newcommand{\sgn}{\mathop{\mathrm{sgn}}}

\begin{document}

\title{Parallel Higher-order Boundary Integral Electrostatics Computation on Molecular Surfaces with Curved Triangulation}

\author{Weihua Geng 
\footnote{ Tel: (205) 348-5302, Fax:
(205)-348-7067, Email: wgeng@as.ua.edu}
\\
\\
\small \it     Department of Mathematics, University of Alabama, Tuscaloosa, AL 35406, USA \\
}

\maketitle

\begin{abstract}
In this paper, we present a parallel higher-order boundary integral method to solve the linear Poisson-Boltzmann (PB)  equation. 
In our method, a well-posed boundary integral formulation is used to ensure the fast convergence of Krylov subspace linear solver such as GMRES.  
The molecular surfaces are first discretized with flat triangles and then converted to curved triangles with the assistance of normal information at vertices.   
To maintain the desired accuracy, four-point Gauss-Radau quadratures are used on regular triangles and sixteen-point Gauss-Legendre quadratures together with regularization transformations  are applied on singular triangles. To speed up our method, we take advantage of the embarrassingly parallel feature of boundary integral formulation, and parallelize the schemes with the message passing interface  (MPI) implementation.  Numerical tests show significantly improved accuracy and convergence of the proposed higher-order boundary integral Poisson-Boltzmann (HOBI-PB) solver compared with boundary integral PB solver using often-seen centroid collocation on flat triangles. The higher-order accuracy results achieved by present method
are important to sensitive solvation analysis of biomolecules, particularly when accurate electrostatic surface potentials are critical in the molecular simulation. In addition, the higher-order boundary integral schemes  presented here and their associated parallelization potentially can be applied to solving boundary integral equations in a general sense. 
\vspace*{1cm}


\end{abstract}

\newpage

\section{Introduction}
Molecular modeling is  
a rising interdisciplinary approach on the study of structure, 
function and dynamics of molecules with biological significance \cite{schlick-02}. 
Among interactions in molecular modeling, 
electrostatics are critical due to their ubiquitous existence.  
Meanwhile, electrostatics are expensive to compute as 
they are long-range and pairwise interactions. 
The Poisson-Boltzmann (PB) model is an effective approach 
to resolve the electrostatics including energy, potential and force
of solvated biomolecules \cite{BakerCurrOp}. 
As an implicit solvent model, 
the PB model considers solvent effects  with a mean field approximation, 
and models the mobile ions with the Boltzmann distribution. 
These implicit treatments of solvent surroundings 
make the PB model computationally more efficient compared with explicit solvent models, 
in which atomic details of solvent molecules and electrolytes are described. 
Recently, experimentalists also showed interests in using the PB model to provide
references for newly released structures of biomolecules e.g. the neurotransmitter receptor Acetylcholine  \cite{Unwin_2003} and the Circadian clock complex CLOCK:BMAL1 \cite{Huang}. 

\begin{figure}[tbh]
  \begin{center}
  \flushleft
  (a)\hs{9.2cm}(b)\hs{1cm}\\
  \includegraphics[width=3.0in]{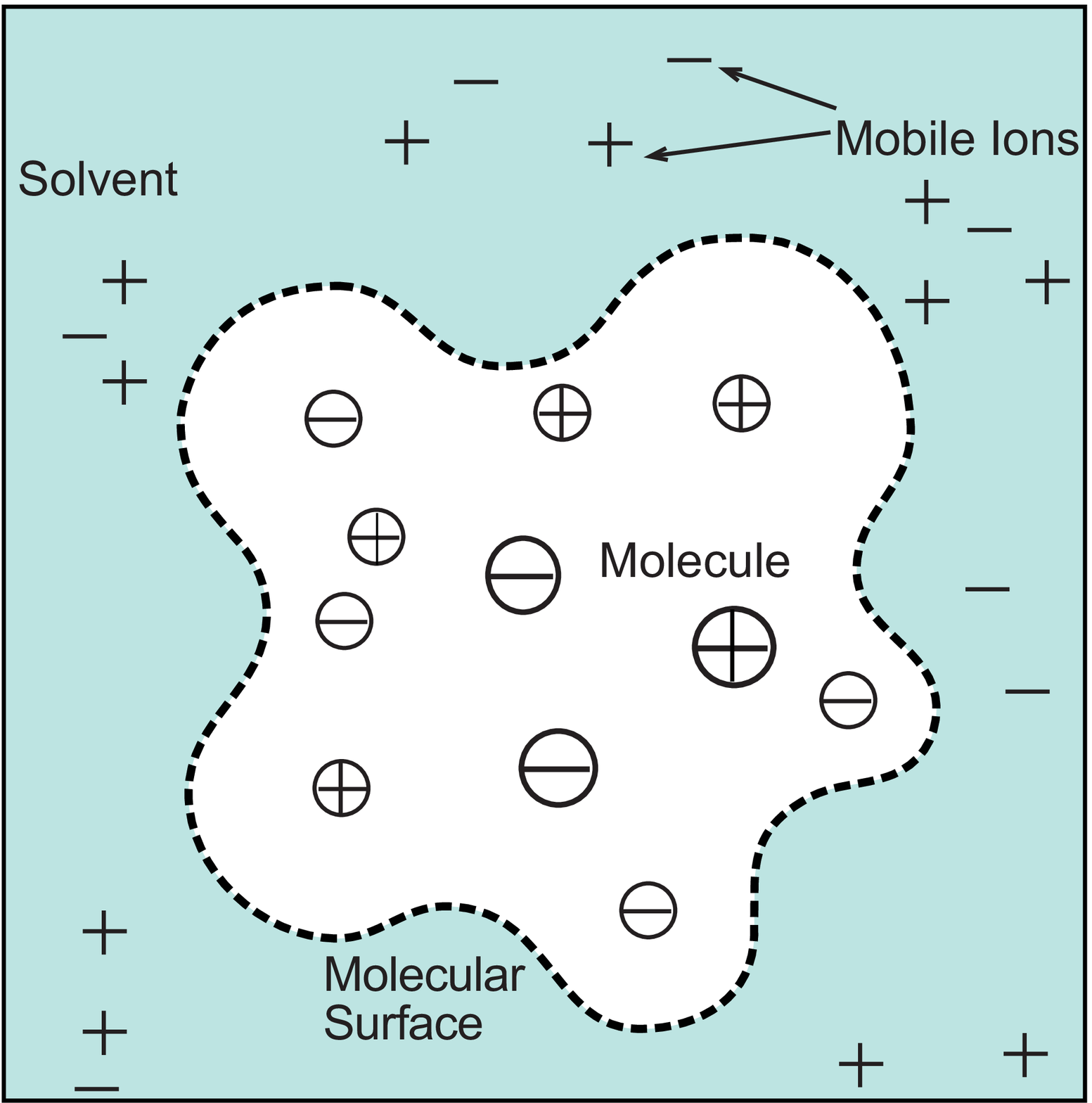} 
  \hs{2cm}
  \includegraphics[width=2.7in]{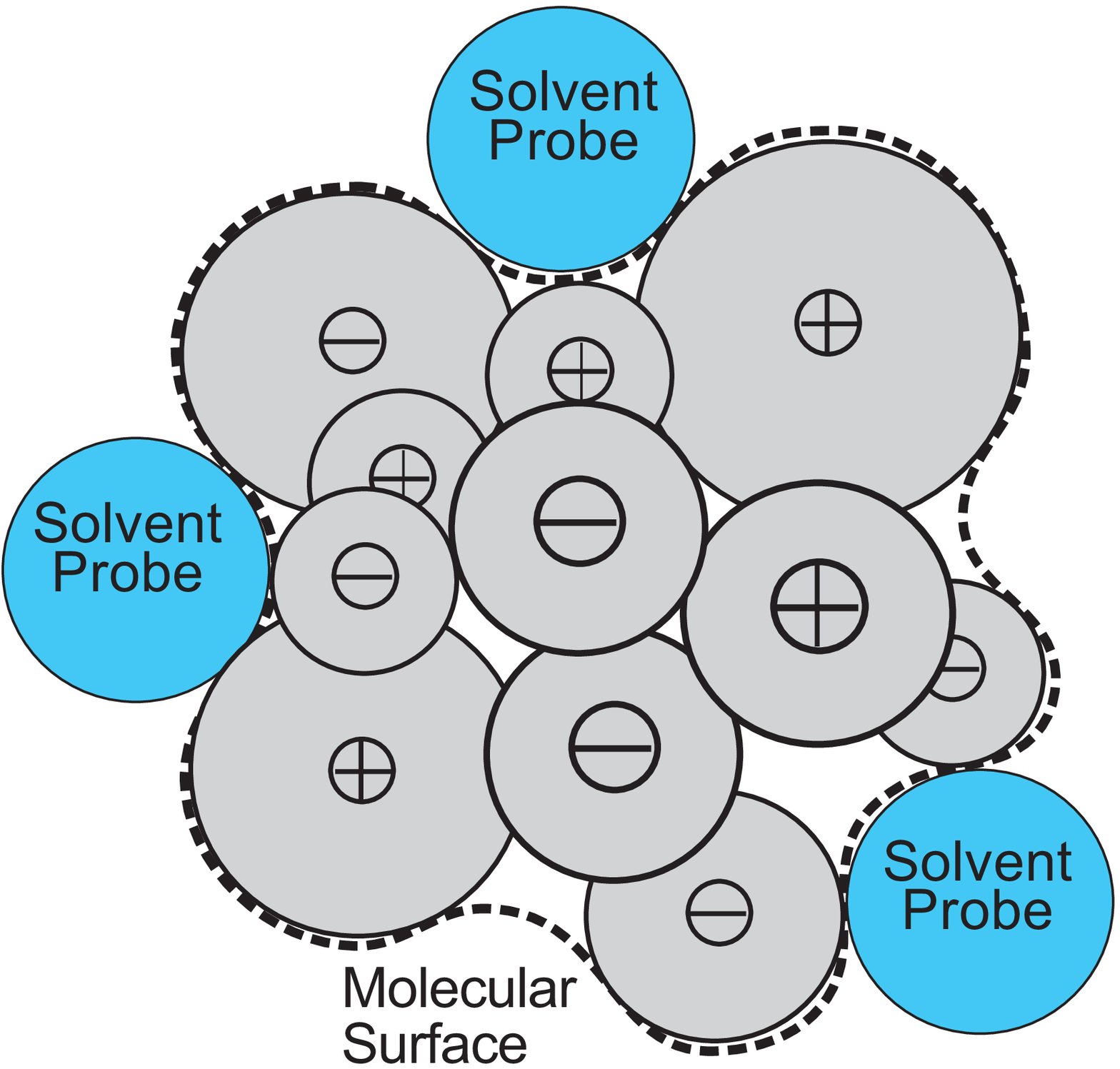}\\
  \end{center}
  \setlength{\unitlength}{.1in}
  \vskip -3.5in
  \begin{picture}(35,35)(0,0)
  \put(0.5,27.0){\makebox{{\huge$\Omega_2$}:{\Large $\varepsilon_2$}}}
  \put(16.5,16.2){\makebox{{\huge$\Omega_1$}:{\Large $\varepsilon_1$}}}
   \put(7,2){\makebox{{\huge$\Gamma$}}}
   \put(45,1.7){\makebox{{\huge$\Gamma$}}}
  \end{picture}
\vskip -10pt
\caption{\small (a) Poisson-Boltzmann (PB) model: domains $\Omega_1 (\text{molecule})$ and $ \Omega_2 (\text{solvent})$ with different dielectric constants $\varepsilon_1$ and $\varepsilon_2$ respectively; (b) the molecular surface is formed by the trace of solvent probe in contact with the solute (molecule). }
\label{model}
\end{figure}

In the PB model, as illustrated in Fig.~\ref{model}(a), the computational domain $\mb{R}^3$ 
is divided into the solute (molecule) domain $\Om_1$ and 
the solvent domain $\Om_2$ by a closed molecular surface $\Gamma$ 
such that $\mb{R}^3=\Om_1\cup\Om_2\cup\Gamma$. 
As shown in Fig.~\ref{model}(b), the molecular surface $\Gamma$ is formed 
by the traces of a spherical solvent probe 
rolling in contact with the van del Walls balls of the solute atoms~\cite{Richards, Connolly85}.  
The solvated molecule, which is located  in domain $\Om_1$ 
with dielectric constant $\veps_1$, 
is represented 
by a set of $N_c$ point charges carrying 
$Q_i$ charge in the units of $e_c$ at positions $\txf{x}_i, i=1, ..., N_c$.
The exterior domain 
contains the solvent with dielectric constant $\veps_2$,
as well as mobile ions.
For $\txf{x}=(x, y, z)$,
the PB equation for the electrostatic potentials in each domain
is derived from the Boltzmann distribution and Gauss' law
and has the form
\begin{align}
\label{eq1}
&\na\cdot(\veps_1(\txf{x})\na\ph_1(\txf{x}))=-\sum_{i=1}^{N_c}
q_i\de (\txf{x}-\txf{x}_i)\q\text{in} ~\,\Om_1,\\
\label{eq2}
&\na\cdot(\veps_2(\txf{x})\na\ph_2(\txf{x}))-\ka^2\sinh\ph_2(\txf{x})=0\q
\text{in} ~\,\Om_2,\\
\label{jump}
&\ph_1(\txf{x})=\ph_2(\txf{x}),\q \veps_1
\fr{\pa\ph_1(\txf{x})}{\pa\nu}=\veps_2\fr{\pa\ph_2(\txf{x})}{\pa\nu}\q
\text{on} ~\,\Ga,\\
\label{radiation}
&\lim_{|\txf{x}|\to\iy}\ph_2(\txf{x})=0,
\end{align}
where $\ph_1$ and $\ph_2$ are the electrostatic potentials in each domain,
$q_i=e_c Q_i/k_B T, i=1, ..., N_c,$
$e_c$ the electron charge,
$k_B$ the Boltzmann's constant,
$T$ the absolute temperature,
$\delta$ the Dirac delta function,
$\ka$ the Debye--H\"{u}ckel parameter, and
$\nu$ the unit outward normal on the interface $\Ga$.
We assume weak ionic strength in this context therefore the non-linear $\sinh$ function term can be approximated by its linearized term, resulting in the linear PB equation 
with $\sinh \phi_2(\bf x)$ term replaced by $\phi_2(\bf x)$ in Eq.~(\ref{eq2}).

The linear PB model is an elliptic equation defined on multiple domains
with discontinuous coefficients across the domain interfaces.
The PB equation has an analytical solution only for the simple geometries
such as spheres \cite{Kirkwood34} or rods \cite{Holst}.
For molecules with complex geometries, the PB equation can only be solved numerically,
which is challenging due to the following numerical difficulties. \\
(1) The solutions to the PB equation, physically the electrostatic potentials,
are not smooth across the interface
as the continuities of both the potentials and the fluxes in Eq.~(\ref{jump}) across the interface $\Gamma$ are required to be satisfied. \\
(2) The complex geometry of the interface needs to be captured
to maintain the accuracy of the potentials particularly on or near the interface. \\
(3) The partial charges carried by the individual atoms of the solute,
modeled by the weighted summation of the Dirac delta functions, is hard to accurately discretize.\\
(4) The PB equation is defined on the entire $\mathbb{R}^3$ domain
subject to boundary condition that the potentials approach zero at infinity thus
a cutoff for 3D mesh-based methods is inevitable. 

The wide application of the PB model as well as its associated numerical difficulties
attracted attention from various computational science communities ranging from
biophysics, biochemistry, mathematics, computer science, mechanical engineering as well as electrical engineering.
Many numerical PB solvers were developed  and they can be roughly but not completely divided into two categories:
The 3D mesh-based finite difference/finite element methods
\cite{Rocchia01,Im98,Luo02,bakerSept,Qiao,Chenjcc11};
and the boundary integral methods
\cite{JBKPB, BFZ, Lu06, yokota-bardhan-knepley-barba-hamada-11,altman-bardhan-white-tidor-09,GengKrasnyjcp12, bajaj-chen-rand-11, Bordner, LS, VorSch}.
All these methods have their own advantages and disadvantages.
For example, the PB solvers embedded in molecular modeling packages
such as Dephi \cite{Rocchia01}, CHARMM \cite{Im98}, AMBER \cite{Luo02}, APBS \cite{bakerSept}
use standard seven-point finite difference with approximated approaches 
to bypass the numerical difficulties (1)-(4). 
Although arguably these solvers have reduced accuracy,
the efficient, robust and user-friendly features of these PB solvers
brought their popularities among the bio-oriented community. 
An often overlooked drawback of these solvers is that 
they provide acceptable accuracy of resolved electrostatics potentials away from the interface 
but are unable to provide accurate solutions \underbar{near} or \underbar{on} the interfaces. 
Applications such as molecular simulation related to ion channels \cite{baker01} 
, cell membranes \cite{FKLS}, and chromatin packing\cite{beard-schlick-01a} 
require accurate electrostatic potentials and fields near or on the interface, 
thus call for developing higher-order methods to solve the PB equation.
Three-dimensional mesh-based Interface methods such as Immersed Interface Methods (IIM) \cite{Qiao}
and Matched Interface and Boundary Poisson-Boltzmann (MIBPB) solver \cite{Chenjcc11}
can significantly improve the accuracy by rigorously treating the numerical difficulties (1)-(3).
However, these methods need to cope with numerical difficulty (4) and the complexities of the algorithms often reduce the efficiency.

Compared with 3D mesh-based methods, the boundary integral methods have many
advantages.\\ 
(1) The solution is characterized solely in terms of
surface distributions so there are fewer unknowns in comparison
to methods that discretize the entire domain. \\
(2) The far-field boundary condition in Eq.~(\ref{radiation}) is exactly imposed. \\
(3) The surface geometry of $\Gamma$ can be represented to high precision using appropriate boundary elements.\\
(4) The electrostatic potential at charge sites are accurately determined using exact
analytical expressions.\\
(5) The continuity conditions in Eq.~(\ref{jump}) are explicitly enforced.\\
Due to these advantages,
boundary integral PB solvers have gained increased attention and we briefly review studies relating to the present work.
In 1988, Zauhar et al.~\cite{ZauMor1988} introduced the boundary integral formulation
by solving the Poisson equation for the induced surface charges. In 1990, Yoon et al.~\cite{Yoon} formulated an ill-posed integral formulation of the PB equation. 
In 1991, Juffer et al.~\cite{JBKPB} reformulated the previous work to obtain a well-posed formulation, which were applied by most of the boundary integral PB solvers after that. 
The boundary integral methods can analytically circumvent the numerical difficulties (1)-(4),
and accelerate the solver with fast algorithms such as fast multipole method (FMM) 
\cite{BFZ,Lu06,yokota-bardhan-knepley-barba-hamada-11}
and treecode \cite{GengKrasnyjcp12}. These boundary integral PB solvers mostly applied centroid collocation methods on flat triangle and benchmark tests on spherical cavities with available analytical solutions show 0.5th order accuracy relative to number of elements \cite{Lu06, GengKrasnyjcp12}, which left spaces for the more challenging problem of developing higher-order boundary integral PB solver.  

In this paper, we present a more accurate boundary integral PB solver on curved triangles with higher-order quadratures and regularization of singularities. The rest of the paper is organized as follows. In section 2, we provide our algorithms including the well-posed boundary integral formulation and the higher-order numerical schemes, followed by the MPI parallelization. In section 3, we provide the numerical results, first on the spherical cavities with centered and eccentric partial changes whose analytical solutions are available and then  on a protein (PDB: 1ajj) for the electrostatics solvation energy computation and the parallel efficiency. This paper ends with a section of concluding remarks.   

\section{Methods}

We will use the well-posed boundary integral formulation 
from Juffer's work \cite{JBKPB} together 
with high order quadrature \cite{ZauMor1988} 
and singularities regularization by a transformation \cite{AtkinsonUserGuide, Schwab}.
We modify and improve these methods as needed and we will explain the details in this section. 
One fact affecting the accuracy of boundary integral methods is the discretization of the surface. 
For the sphere, we use a non-uniformed triangular surface from MSMS \cite{Sanner} 
with radial projection to correct the truncated output,  or a uniform icosahedral triangulation  \cite{Baumgardner}. For biomolecules, we only use MSMS to generate the flat triangles. All flat triangles are then converted to curved triangles by applying the schemes introduced as follows.  Through the paper, we call our higher-order boundary integral Poisson-Boltzmann solver as HOBI-PB solver. 

\subsection{Well-posed integral formulation}
The differential PB equation in Eqs.~(\ref{eq1}) and~(\ref{eq2}) can be converted to boundary integral equation. By applying the fundamental solution of Poisson equation~(\ref{eq1}) , $G_0$, in $\Omega_1$ and the fundamental solution of PB equation~(\ref{eq2}), $G_{\kappa}$, in $\Omega_2$, together with Green's second theorem, and cancel the normal derivative terms with interface jump conditions in Eq.~(\ref{jump}), the coupled integral equations can be derived as \cite{Yoon}:
\begin{align}
\label{eqbim_1}
\ph_1(\txf{x})=&\int_{\Ga} \left[G_0(\txf{x},
\txf{y})\fr{\pa\ph_1(\txf{y})}{\pa\nu_{\txf{y}}} -\fr{\pa
G_0(\txf{x}, \txf{y})}{\pa\nu_{\txf{y}}} \ph_1(\txf{y})
\right]dS_{\txf{y}}+ \sum_{k=1}^{N_{c}}
q_k G_0(\txf{x}, \txf{y}_k) ,\q\hskip 0.1in\txf{x}\in\Om_1,\\
\label{eqbim_2}
\ph_2(\txf{x})=&\int_{\Ga}\left[-G_\kappa(\txf{x},
\txf{y})\fr{\pa\ph_2(\txf{y})}{\pa\nu_{\txf{y}}}+\fr{\pa
G_\kappa(\txf{x}, \txf{y})}{\pa\nu_{\txf{y}}}
\ph_2(\txf{y})\right]dS_{\txf{y}}, \q\hskip 1.18in \txf{x}\in\Om_2.
\end{align}
where $G_0(\txf{x},\txf{y})$ and $G_\kappa(\txf{x}, \txf{y})$ are the 
Coulomb and screened Coulomb potentials,
\begin{equation}
G_0(\txf{x}, \txf{y}) = \fr{1}{4\pi|\txf{x}-\txf{y}|},
\quad
\label{eq_potential}
G_\kappa(\txf{x}, \txf{y}) = \fr{e^{-\kappa|\txf{x}-\txf{y}|}}{4\pi|\txf{x}-\txf{y}|}.
\end{equation}
However,
straightforward discretization of Eqs.~(\ref{eqbim_1}) and~(\ref{eqbim_2}) yields a linear system which becomes
ill-conditioned as the number of boundary elements increases~\cite{Lu07}. Juffer et al. derived a well-posed boundary integral formulation by going through the differentiation of the single-layer and double-layer potentials \cite{JBKPB}. The desired forms are:
\begin{align}
\label{eqbim_3}
\fr{1}{2}\left(1+\veps\right)\ph_1(\txf{x})&=
\int_{\Ga}\left[K_1(\txf{x}, \txf{y})\fr{\pa\ph_1(\txf{y})}
{\pa\nu_{\txf{y}}}+K_2(\txf{x}, \txf{y})\ph_1(\txf{y})\right]dS_{\txf{y}}+S_{1}(\txf{x}),
\qquad \txf{x}\in\Ga, \\
\label{eqbim_4}
\fr{1}{2}\left(1+\fr{1}{\veps}\right)\fr{\pa\ph_1(\txf{x})}{\pa\nu_{\txf{x}}}&=
\int_{\Ga}\left[K_3(\txf{x}, \txf{y})\fr{\pa\ph_1(\txf{y})}
{\pa\nu_{\txf{y}}}+K_4(\txf{x}, \txf{y})\ph_1(\txf{y})\right]dS_{\txf{y}}
+S_{2}(\txf{x}),
\qquad \txf{x}\in\Ga,
\end{align}
with the notation
\begin{align}\nonumber
K_1(\txf{x}, \txf{y})=&\,{G_{0}(\txf{x},\txf{y})}-{G_{\kappa}(\txf{x},\txf{y})}, \hskip 0.75in
K_2(\txf{x}, \txf{y})=\veps\fr{\pa G_{\kappa}(\txf{x},\txf{y})}{\pa\nu_{\txf{y}}}-\fr{\pa
G_{0}(\txf{x},\txf{y})}{\pa\nu_{\txf{y}}},\\
\label{Eq_Ls}
K_3(\txf{x}, \txf{y})=&\,\fr{\pa
G_{0}(\txf{x},\txf{y})}{\pa\nu_{\txf{x}}}-\fr{1}{\veps}\fr{\pa
G_{\kappa}(\txf{x},\txf{y})}{\pa\nu_{\txf{x}}}, \hskip 0.4in
K_4(\txf{x}, \txf{y})=\fr{\pa^2
G_{\ka}(\txf{x},\txf{y})}{\pa\nu_{\txf{x}}\pa\nu_{\txf{y}}}-\fr{\pa^2
G_{0}(\txf{x},\txf{y})}{\pa\nu_{\txf{x}}\pa\nu_{\txf{y}}},\\
S_{1}(\txf{x})=&\,\sum_{k=1}^{N_{c}}q_kG_{0}(\txf{x}, \txf{y}_k), \hskip 1.2in
S_{2}(\txf{x})=\sum_{k=1}^{N_{c}}q_k
\fr{\pa
G_{0}(\txf{x},\txf{y}_k)}{\pa\nu_{\txf{x}}}\label{Eq_Ss}.
\end{align}
and $\veps={\veps_{1}}/{\veps_{2}}$.
Note this is the well-posed Fredholm second kind of integral equation which is also our choice in this paper.

In order to numerically solve the coupled equations~(\ref{eqbim_3}) and (\ref{eqbim_4}), we need to discretize the molecular surface $\Gamma$ with high quality elements and implement the numerical integral with higher-order quadrature. We also need to treat the occurred singularities or near-singularities when $\bf x$
and $\bf y$ are equal or nearly equal in Kernels $K_{1,\dots, 4}$. These details are 
described in the following subsections. 
\subsection{Curved triangles and higher-order quadratures}

\begin{figure}
\centering
\includegraphics[width=1.8in]{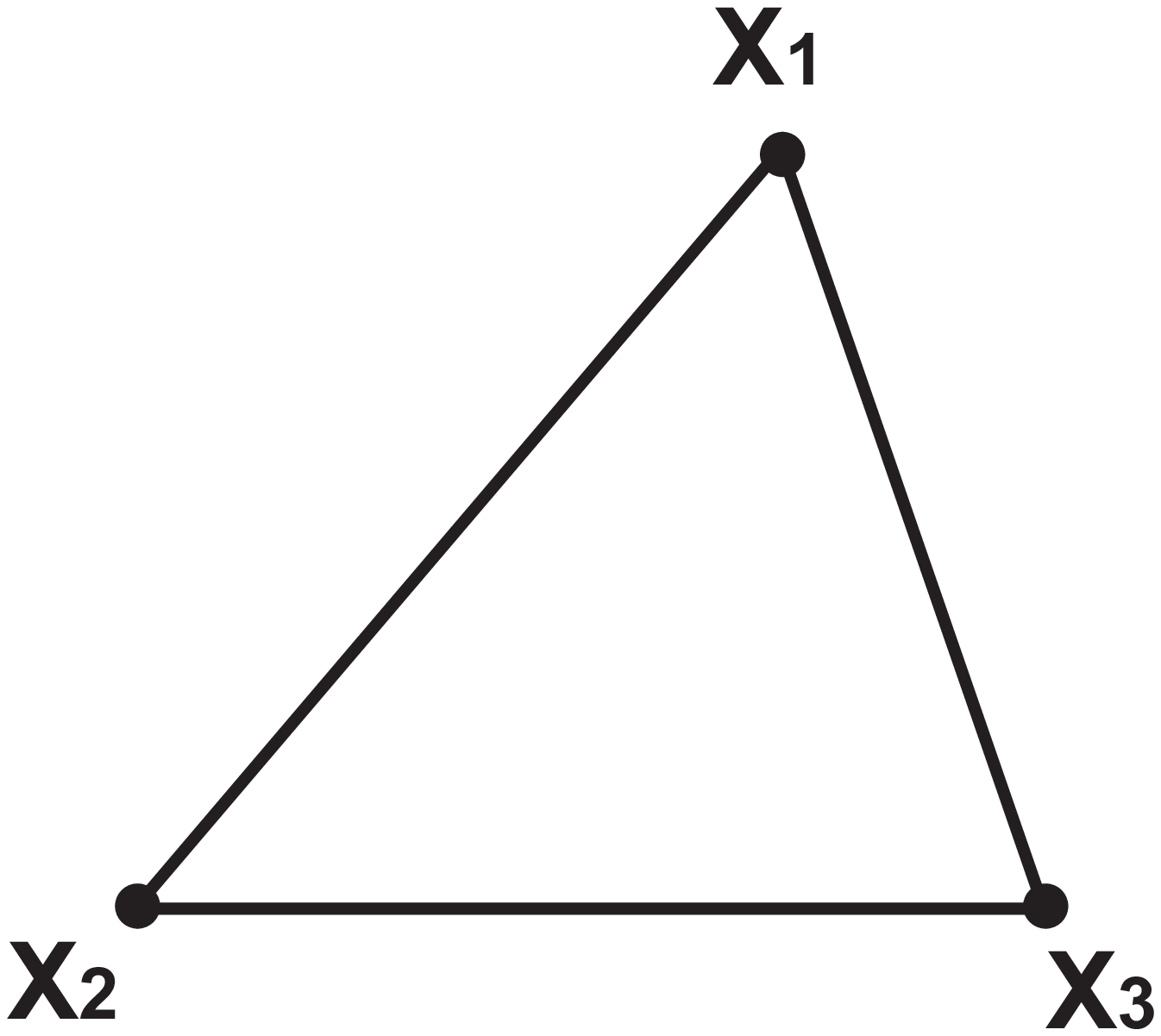}
\hs{1cm}
\includegraphics[width=1.8in]{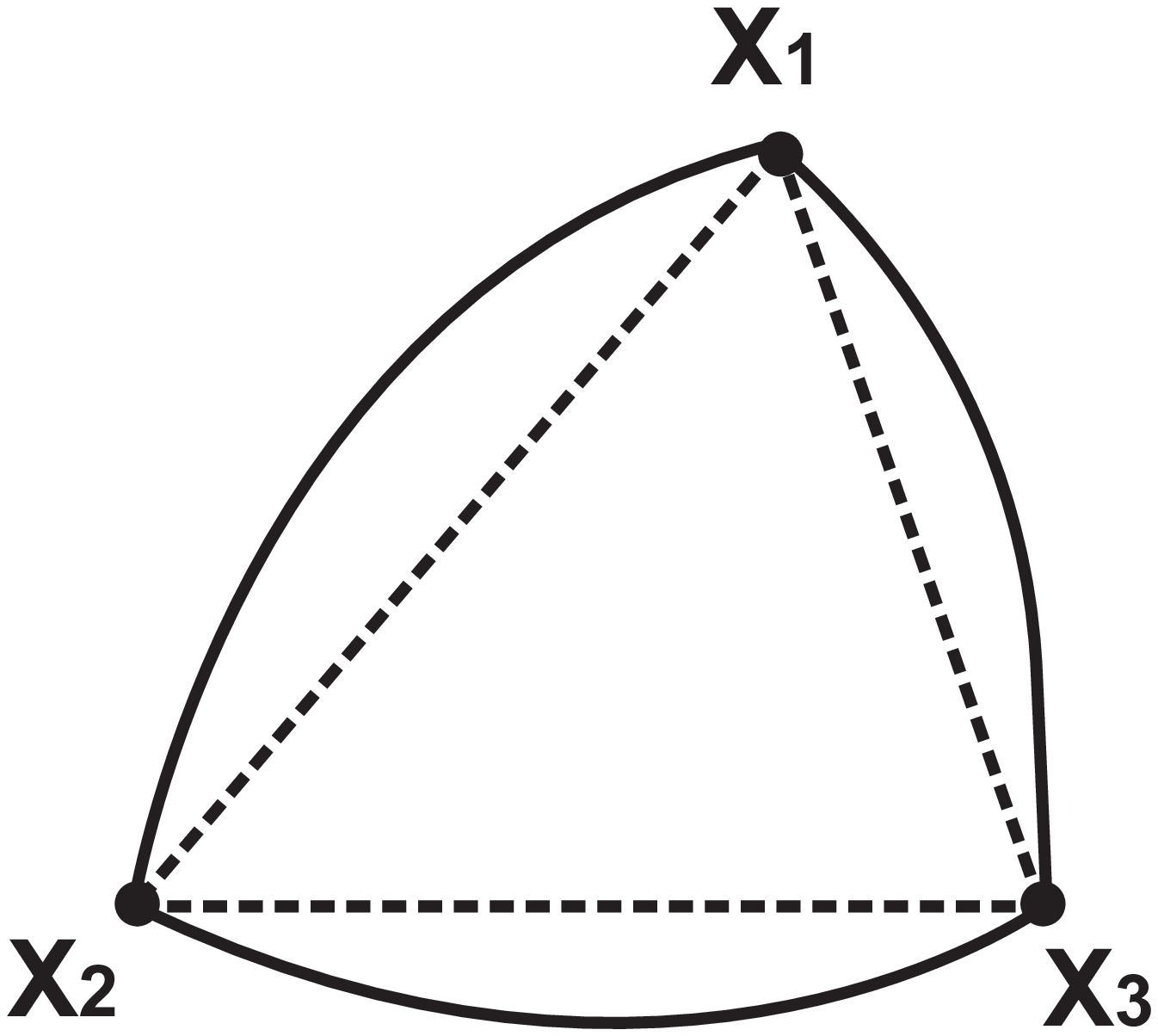}\hs{1cm}\includegraphics[width=1.8in]{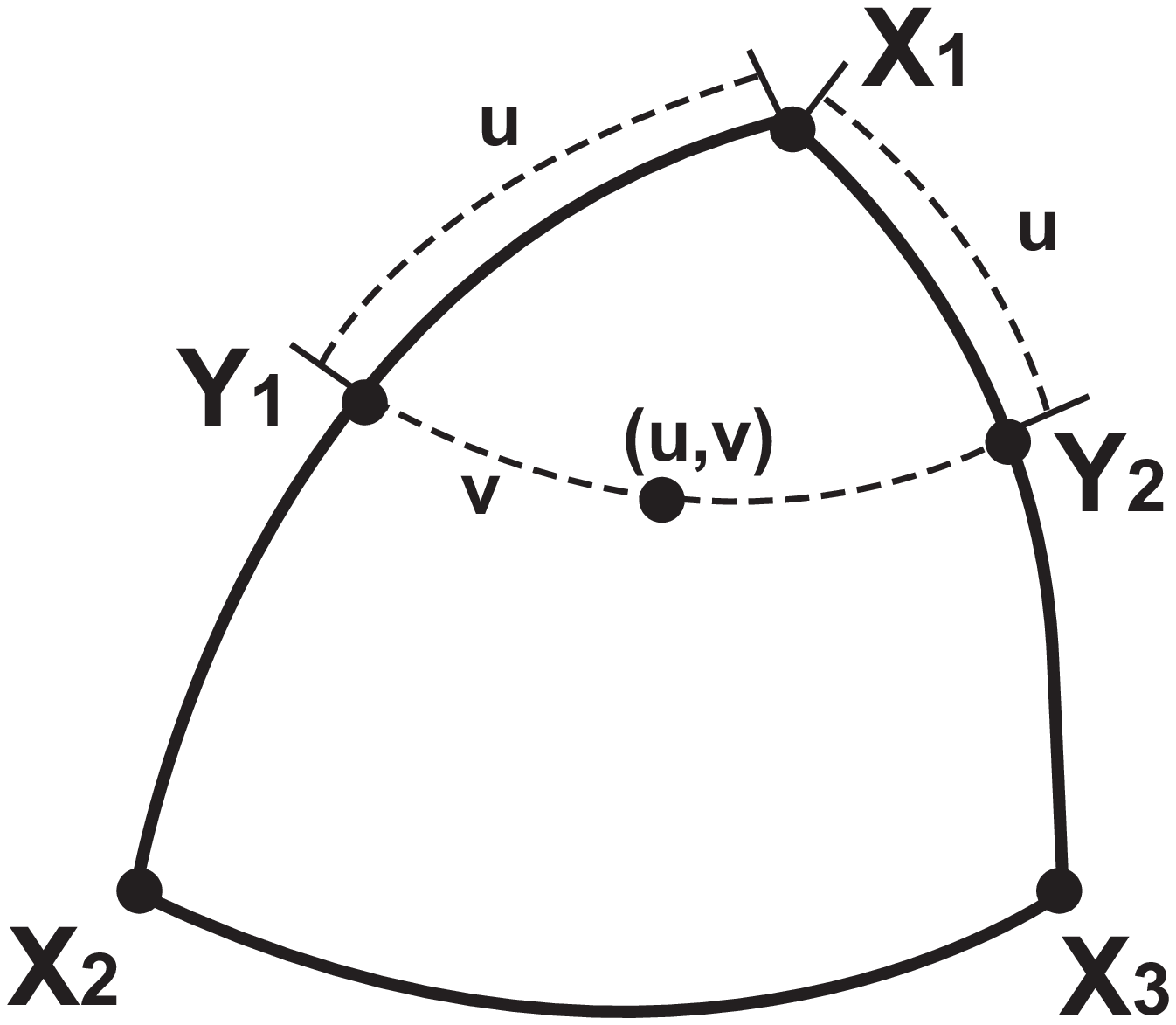}\\
(a)\hs{6cm}(b)\hs{6cm}(c)\\
\vskip -10pt
\caption{Triangular surface elements: (a) a flat triangle (b)a curved triangle (c) curve linear coordinates (u,v) in the curved triangle.}\label{fig_triangle}
\end{figure}

Given the positions and radii of all atoms of a molecule, a triangulation  program 
e.g. MSMS \cite{Sanner} can generate a discretized surface with a set of $N_f$ flat triangles, $N_v$ nodes (vertices) and corresponding normal directions. 
To achieve higher-order accuracy, we apply the schemes described in \cite{ZauMor1988} 
to convert these flat triangles to curved triangles 
under the new background of solving the well-posed integral PB equation in Eqs.~(\ref{eqbim_3}) and (\ref{eqbim_4}). 
We also modify and improve these schemes 
to treat the singularities of kernels $K_{1,\dots, 4}$ in Eq.~(\ref{Eq_Ls}). 
To keep this paper in an integrated form, 
we start from restating the schemes in \cite{ZauMor1988} to produce the curved triangles and higher-order quadratures for discretizing Eqs.~(\ref{eqbim_3}) and (\ref{eqbim_4}).  

We first replace the straight element edges as shown in Fig.~\ref{fig_triangle}(a) with curved arcs in Fig.~\ref{fig_triangle}(b) \cite{ZauMor1988}.
Let $\txf{x}(t)$ be the arc between two nodes,
parameterized by the dimensionless variable $t$,
\begin{equation}
\txf{x}(t)=\txf{c}_0+\txf{c}_1 t+\txf{c}_2 t^2+\txf{c}_3 t^3,
\end{equation}
where $\txf{c}_0, \txf{c}_1, \txf{c}_2, \txf{c}_3$ are vector constants (12 unknowns), which will be
determined by a pair of connected nodes and associated unit normals (12 conditions). 
At any point on the curve, 
the normal can be found by ${\bf n}(t)=\sgn(t)\displaystyle\frac{\kappa(t)}{|(\kappa(t))|}$, 
where $\sgn(t)=\pm1$ is chosen to keep a constant orientation of ${\bf n}(t)$ along the curve 
and $\kappa(t)$ is the curvature given as \cite{ZauMor1988}:
\begin{equation}\label{eqcurvature}
    \kappa(t)=\displaystyle\frac{1}{d{\bf x}/dt}\left[\frac{d^2 {\bf x}}{dt^2}-\frac{(\frac{d{\bf x}}{dt}\cdot\frac{d^2{\bf x}}{dt^2})}{|d{\bf x}/dt|^2}\frac{d{\bf x}}{dt}\right]
\end{equation}

As shown in Fig.~\ref{fig_triangle}(c), we can use parameter $u\in [0,1]$ for the two curves starting from $X_1$ and ending at $X_2$ and $X_3$. Then for any given $u$, two points on the curves $X_1X_2$ and $X_1X_3$ are specified, say $Y_1$ and $Y_2$ with normal directions.  By following the same procedure, we could find a curve connecting $Y_1$ and {$Y_2$}, using $v\in [0,1]$ as the parameter. The pair of parameters $(u,v)$ will therefore correspond to a point on the curve element.

\begin{figure}
  \begin{center}
  \includegraphics[width=6.5in]{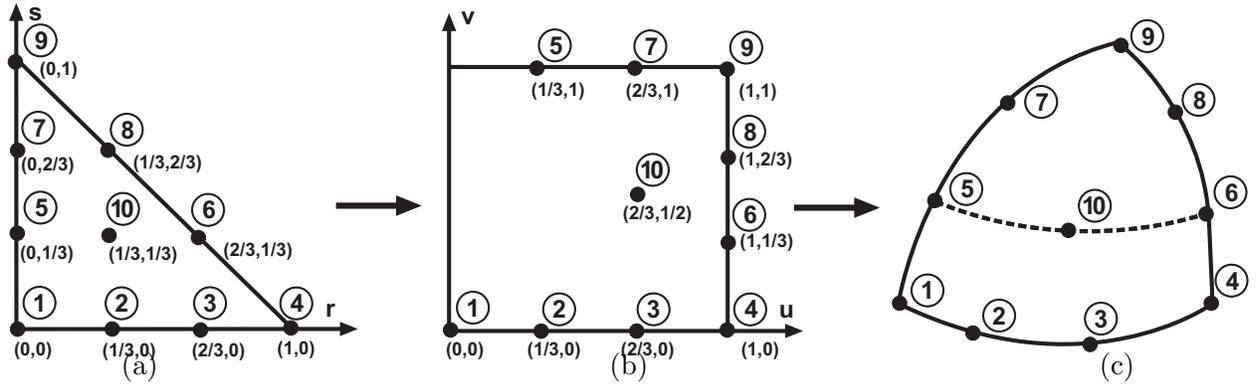}\\
\vskip -10pt
  (a)\hs{6cm}(b)\hs{6cm}(c)\\
  \caption{Coordinates Transformation}\label{fig_all10pts}
  \end{center}
\end{figure}

In order to conveniently use quadrature rules, the integral is conducted in a unit right triangle. A mapping between the parameters $r,s \in [0,1]$ on the unit right triangle in Fig.~\ref{fig_all10pts}(a) and parameters $u,v$ on a unit square in Fig.~\ref{fig_all10pts}(b) is constructed with the following transformation \cite{ZauMor1988}.
\begin{equation}\label{eqrs2uv}
\left\{
\begin{array}{lr}
u=(r+s), v=s/(r+s) &\text{if } r+s\neq 0\\
u=0,v=0 &\text{if } r=s=0
\end{array}\right.
\end{equation}
Based on this, we can establish the mapping from point $(r,s)$ 
on the unit right triangle in Fig.~\ref{fig_all10pts}(a) to 
point ${\bf x}^{(j)}(u(r,s),v(r,s))={\bf x}^{(j)}(r,s)$ on the $j$th curved elements  in Fig.~\ref{fig_all10pts}(c) by the following steps:\\
(1) Given $(r,s)$, find $u,v$ through Eq.~(\ref{eqrs2uv}).\\
(2) Plug $u$ as the parameter into the curve functions of $X_1X_2$ and $X_1X_3$ to locate $Y_1$ and $Y_2$ as shown in Fig.~\ref{fig_triangle}(c).\\
(3) Find the curve function of $Y_1Y_2$ and then plug $v$ as the parameter to finally get ${\bf x}^{(j)}$ on the curved elements.\\
However, this is not an analytical function for efficient computations but a constructive procedure. As a remedy, a high-order 10-point interpolation scheme is used \cite{ZauMor1988}. The brief idea is illustrated in Fig.~(\ref{fig_all10pts}) on the $j$th element of the triangulation.\\
(1) Pick 10 specified points $(r_k,s_k)$ first in the right unit triangle for $k=1,2,...,10$.\\
(2) Find $u(r_k,s_k),v(r_k,s_k)$ according to Eq.~(\ref{eqrs2uv}).\\
(3) Find ${\bf x}^j(u(r_k,s_k),v(r_k,s_k))$ and their normal directions on the curved element, which is parameterized by $(u,v)$. Note points 1,4,9 are already given with the flat triangles. With these three nodes, three trajectories can be found and used to find the positions and normal directions of points 2,3,5,6,7,8. Finally with points 5 and 6, the dashed trajectory can be formed to find point 10 and its normal direction.\\
(4) As required by the quadrature rules, interpolate any point on the curved element ${\bf x}^{(j)}(r,s)$ by the expression
\begin{equation}\label{eqintepolation}
    {\bf x}^{(j)}(r,s)=\sum\limits_{k=1}^{10}N_k(r,s){\bf x}_k^{(j)}
\end{equation}
where $N_k(r,s)$ are Lagrangian interpolation polynomial for a 10-point element. See table 1 of \cite{ZauMor1988} for the expression of $N_k(r,s)$. Note Eq.~(\ref{eqintepolation}) can also be used to find the partial derivative of ${\bf x}^{(j)}$ with respect to $r$ and $s$, which is required for computing the Jacobian for transformation and the normal direction at ${\bf x}^{(j)}$. 

Suppose now we will integrate function $f({\bf x})$ on the $j$th curved element. The quadrature rules give the position of a set of points on the unit right triangle e.g. $(r_m,s_m)$ and quadrature weights $W_m$ for $m=1,2,..,N$, where $N$ is the number of quadrature points, and the integral can be evaluated as:
\begin{eqnarray}
  \int_{\triangle_j} f({\bf x})dA   
                        &=& \sum\limits_{m=1}^N f({\bf x}(r_m,s_m))|\displaystyle\frac{\partial{\bf x}(r_m,s_m)}{\partial r}\times\frac{\partial{\bf x}(r_m,s_m)}{\partial s}|W_m
\end{eqnarray}
Note the term $\frac{\partial{\bf x}(r_m,s_m)}{\partial r}\times\frac{\partial{\bf x}(r_m,s_m)}{\partial s}$ gives the normal direction of the point ${\bf x}^{(j)}(r,s)$ with parameters $(r,s)$, and we will use this information to supply 
the required normal direction in Eqs.~(\ref{eqbim_3}) and~(\ref{eqbim_4}).

In this paper, we use $N=N_{GR}$ points (practically we choose $N_{GR}=4$) in the Gauss-Radau quadrature \cite{Akin}. For $i=1,2,...,N_v$, the $i$th and the $(i+N_v)$th element of the discretized matrix-vector product ${\bf Au}$ are given as

{\footnotesize
\begin{eqnarray}\label{eqdbef3}
  \{{\bf Au}\}_i^r &=& \fr{1}{2}\left(1+\veps\right)\ph_1(\txf{x}_{i})-\sum\limits_{j=1}^{N_i^r}\sum\limits_{m=1}^{N_{GR}}\sum\limits_{n=1}^3 W_{j,m,n}\left[K_1(\txf{x}_{i}, \txf{y}_{j,m})\fr{\pa\ph_1(\txf{y}_{j,n})}{\pa\nu_{\txf{y}_{j,n}}}+K_2(\txf{x}_{i},\txf{y}_{j,m})\ph_1(\txf{y}_{j,n})\right] \\
   \{{\bf Au}\}_{i+N_v}^r&=&\fr{1}{2}\left(1+\fr{1}{\veps}\right)\fr{\pa\ph_1(\txf{x}_{i})}{\pa\nu_{\txf{x}_{i}}}-\sum\limits_{j=1}^{N_i^r}\sum\limits_{m=1}^{N_{GR}}\sum\limits_{n=1}^3 W_{j,m,n}\left[K_3(\txf{x}_{i}, \txf{y}_{j,m})\fr{\pa\ph_1(\txf{y}_{j,n})}
{\pa\nu_{\txf{y}_{j,n}}}+K_4(\txf{x}_{i}, \txf{y}_{j,m})\ph_1(\txf{y}_{j,n})\right]\label{eqdbef4}
\end{eqnarray}
}\\
In Eqs.~(\ref{eqdbef3}) and~(\ref{eqdbef4}), the superscript $r$ of $\{\bf Au\}$ stands for regular triangle, $W_{j,m,n}$ contains the weights and coefficients associated with the quadrature, the transformation Jacobian, and the interpolation coefficients in Eq.~(\ref{eqintepolation}). Note that  $N_i^r$ is the number of regular triangles associated with the $i$th vertex. In these equations, $\txf{y}_{j,n}$ are the same set of nodes as $\txf{x}_{i}$, but $\txf{y}_{j,m}$ are the quadrature points on the $j$th curved element, which are mapped from the predetermined points on the unit right triangle. This mismatch brings difficulty to apply Fast Multipole Method or treecode to accelerate the higher-order scheme. Studies about this issue will be proceeded in our future work.

\subsection{Treatment of singularities}
From Eqs.~(\ref{eqdbef3}) and~(\ref{eqdbef4}), we can see that when node $i$ is one of the vertices of the $j$th curved triangular element, singularity or near singularity occurs at evaluating kernels $K_{1,\dots, 4}$. 
In other words, ${\bf x}_i$ is equal to or nearly equal to ${\bf y}_{j,m}$. 
It can be shown that the singularities in kernels $K_{1,\dots, 4}$ are 
in the order of $\mathcal{O}(\frac{1}{|{\bf x}_i-{\bf y}_{j,m}|})$ \cite{JBKPB}.
To treat these singularities, we use the tensor-product of Gauss-Legendre quadrature on a unit square, 
together with a transformation \cite{AtkinsonUserGuide, Schwab}. 
In this case the mapping from a unit square ($0\leq x,y \leq 1$) 
to a unit right triangle ($0 \leq r,s \leq 1$ and $r+s \leq 1$), 
and to a curved triangle is constructed. 
The mapping that $r=(1-y)x$ and $s=yx$ for $0 \leq x,y \leq 0$ 
is used to remove the singularities in kernels $K_{1,\dots, 4}$
when they appear at 
$(r,s)=(0,0)$ which indicates $x=0$.
To understand this, it can be seen that the Jacobian 
for the transformation from $(r,s)$ to $(x,y)$ is $x$, 
thus it can remove the $\mathcal{O}(\frac{1}{|{\bf x}_i-{\bf y}_{j,m}|})$ type of singularities.   

In the treatment of singularities, if the number of quadrature points used 
in each direction of the Gauass-Lagendre quadrature is $N_{GL}$, 
the $i$th and $(i+N_v)$th element of the discretized matrix-vector product ${\bf Au}$ 
will in addition contain the following singular component

{\footnotesize
\begin{eqnarray}\label{eqdbef5}
\{{\bf Au}\}_{i}^s &=& \fr{1}{2}\left(1+\veps\right)\ph_1(\txf{x}_{i})-\sum\limits_{j=1}^{N_i^s}\sum\limits_{m=1}^{({N_{GL}})^2}\sum\limits_{n=1}^3 W_{j,m,n}\left[K_1(\txf{x}_{i}, \txf{y}_{j,m})\fr{\pa\ph_1(\txf{y}_{j,n})}{\pa\nu_{\txf{y}_{j,n}}}+K_2(\txf{x}_{i},\txf{y}_{j,m})\ph_1(\txf{y}_{j,n})\right] \\
   \{{\bf Au}\}_{i+N_v}^s&=&\fr{1}{2}\left(1+\fr{1}{\veps}\right)\fr{\pa\ph_1(\txf{x}_{i})}{\pa\nu_{\txf{x}_{i}}}-\sum\limits_{j=1}^{N_i^s}\sum\limits_{m=1}^{({N_{GL}})^2}\sum\limits_{n=1}^3 W_{j,m,n}\left[K_3(\txf{x}_{i}, \txf{y}_{j,m})\fr{\pa\ph_1(\txf{y}_{j,n})}
{\pa\nu_{\txf{y}_{j,n}}}+K_4(\txf{x}_{i}, \txf{y}_{j,m})\ph_1(\txf{y}_{j,n})\right]\label{eqdbef6}
\end{eqnarray}
}\\
In Eqs.~(\ref{eqdbef5}) and~(\ref{eqdbef6}), the superscript $s$ of $\{\bf Au\}$ stands for singular triangle, and $W_{j,m,n}$ contains the weights and coefficients associated with the quadrature, the transformation Jacobians (additionally contains the Jacobian of the $r=(1-y)x$ and $s=yx$ mapping), and the interpolation coefficients in Eq.~(\ref{eqintepolation}). 
Note index $j$ has the range up to $N_i^s$, which stands for the number of singular triangles associated with the $i$th vertex. Simulation shows this value is various for different vertices and it can be as big as 15. 
In this paper, we practically choose $N_{GL}=4$ points in each direction to ensure desired accuracy \cite{Akin}. 

\subsection{Low-order scheme}
We also briefly introduce the low-order boundary integral Poisson-Boltzmann (LOBI-PB) solver, 
which is used for comparison in this paper. 
In this low-order scheme, the flat triangle and centroid collocation are used, 
i.e. the quadrature point is located at the center of each triangle. 
This scheme also assumes that the potential and its normal derivative, as well as the kernel functions are uniform on each triangle. When singularity in kernels occurs, the contribution of this triangle in the integral is then simply removed. This scheme is in fact widely used in latest boundary integral methods in solving PB equations to provide convenience on incorporating fast algorithms such as FMM \cite{Lu07} and treecode \cite{GengKrasnyjcp12}. 
For $i=1,2,...,N_f$, the $i$th and the $(i+N_f)$th element 
of the discretized matrix-vector product ${\bf Au}$ are given as
\begin{eqnarray}\label{eqdbef1}
\{{\bf Au}\}_{i} &=& \fr{1}{2}\left(1+\veps\right)\ph_1(\txf{x}_{i})-\sum\limits_{j=1,j\ne i}^{N_f} W_{j}\left[K_1(\txf{x}_{i}, \txf{x}_{j})\fr{\pa\ph_1(\txf{x}_{j})}{\pa\nu_{\txf{x}_{j}}}+K_2(\txf{x}_{i},\txf{x}_{j})\ph_1(\txf{x}_{j})\right] \\
   \{{\bf Au}\}_{i+N_f}&=&\fr{1}{2}\left(1+\fr{1}{\veps}\right)\fr{\pa\ph_1(\txf{x}_{i})}{\pa\nu_{\txf{x}_{i}}}-\sum\limits_{j=1,j\ne i}^{N_f} W_{j}\left[K_3(\txf{x}_{i}, \txf{x}_{j})\fr{\pa\ph_1(\txf{x}_{j})}
{\pa\nu_{\txf{x}_{j}}}+K_4(\txf{x}_{i}, \txf{x}_{j})\ph_1(\txf{x}_{j})\right]\label{eqdbef2}
\end{eqnarray}
It worths mentioning that the unknowns of LOBI-PB is at the centroid of the triangle in the number of $N_f$ while the unknowns of HOBI-PB is at the vertices of the triangle in the number of $N_v$.

\subsection{Electrostatic Solvation Energy Formulation}
The electrostatic solvation energy is computed by
\begin{equation}
E_{\rm sol} = \frac{1}{2}\sum_{k=1}^{N_c}q_k\phi_{\rm reac}(\txf{x}_k) =
\frac{1}{2}\displaystyle \sum\limits_{k=1}^{N_c} q_k
\int_\Gamma \left[K_1(\txf{x}_k, \txf{y})\fr{\pa\ph_1(\txf{y})}
{\pa\nu_\txf{y}}+K_2(\txf{x}_k, \txf{y})\ph_1(\txf{y})\right]dS_{\txf{y}},
\label{solvation_energy}
\end{equation}
where $\phi_{\rm reac}({\bf x}_k) = \phi_1({\bf x}_k)-S_1({\bf x}_k)$, whose formulation is the integral part of Eq.~(\ref{solvation_energy}),  
is the reaction potential at the $k$th solute atom. 
The electrostatic solvation energy,  which can be regarded as the atomistic charge weighted average of the reaction potential $\phi_{\rm reac}$, can effectively characterize the accuracy of a PB solver. 

\subsection{MPI Based Parallel Implementation}

\begin{table}[htb]
\caption{\small 
Pseudocode for parallel HOBI-PB solver using replicated data algorithm.}
\begin{center}
\vskip -7.5pt
\fbox{
\begin{tabular}{rl}
1 & on main processor \\
2 & \qquad read protein data \\
3 & \qquad call MSMS to generate triangulation \\
4 & \qquad copy protein data and triangulation to all other processors \\
5 & on each processor \\
6&  \qquad locally compute sources terms for each assigned vertex\\
7 & \qquad locally convert flat triangles to curved triangles\\
8 & \qquad for each assigned triangle, locally compute and store quadrature information \\
9 & \qquad for each assigned vertex, locally store quadrature of associated singular triangles\\
10 & \qquad\qquad copy result to all other processors \\
11 & \qquad set initial guess for GMRES iteration \\
12 & \qquad compute assigned segment of matrix-vector \\
13 & \qquad\qquad copy result to all other processors \\
14 & on main processor \\
15 & \qquad test for GMRES convergence \\
16 & \qquad\qquad if no, go to step 12 for next iteration \\
17 & \qquad\qquad if yes, go to step 17 \\
18 & on each processor \\
19 & \qquad compute assigned segment of electrostatic solvation energy \\
20 & \qquad\qquad copy result to main processor \\
21 & on main processor \\
22 & \qquad add segments of electrostatic solvation energy and output result \\
\end{tabular}
}\label{tb_MPI}
\end{center}
\label{pseudocode}
\end{table}

The HIBO-PB solver uses the boundary integral formulation, which can be conveniently parallelized. 
The majority of the CPU time is taken by the following routines.\\
(1) Compute the source term in Eqs.~(\ref{eqbim_3}) and~(\ref{eqbim_4}). \\
(2) Convert flat triangles to curved triangles. \\
(3) Compute and store the Gauss-Radau quadrature information for each curved triangle.\\
(4) Compute and store the Gauss-Legendre quadrature information for singular triangles associated with each vertex. \\
(5) Perform matrix-vector product as in Eqs.~(\ref{eqdbef3}), (\ref{eqdbef4}), (\ref{eqdbef5}), and (\ref{eqdbef6}) on each GMRES iteration. \\
(6) Compute electrostatic solvation energy with Gauss-Radau quadrature. \\
Among all these routines, routine (5) is the most expensive one and it will be repeatedly computed on each GMRES iteration. We parallelize all these routines to maximize the parallel efficiency.  Table \ref{tb_MPI}  provides the Pseudocode for MPI based parallel implementation. 

\section{Results}
In this section, we present numerical results. 
We first solve the PB equation 
on spherical cavities with a centered charge 
and with an eccentric change at different locations. 
The analytical solutions in terms of a closed form (for centered charge) and 
in terms of spherical harmonics (for multiple eccentric charges)
are available for these tests \cite{Kirkwood34}. 
To demonstrate the higher accuracy obtained by higher-order boundary integral Poission--Boltzmann (HOBI-PB) solver, 
we compare the numerical results with APBS\cite{bakerSept}, MIBPB\cite{Chenjcc11}, 
and the lower-order boundary integral Poisson-Boltzmann (LOBI-PB) solver. APBS uses straight-forward finite difference scheme. 
MIBPB is a 2nd order interface method repeatedly using local interpolation to capture interface jump conditions  \cite{Zhoujcc08,YuGengjcp07}
and applying a Dirichlet-to-Neumann mapping 
to transform the singular charges to interface jump conditions \cite{Gengjcp07}.  LOBI-PB discretizes the integral equations with flat triangles and performs the numerical integral with centroid collocation as explained in the previous section. 

We then solve the PB equation on a test protein (PDB: 1ajj), 
which is a lipprotein receptor with 37 residues and 519 atoms. 
We report the electrostatic solvation energy results for this protein 
computed from both HOBI-PB and LOBI-PB 
to demonstrate the improved accuracy achieved by the higher-order integral schemes.
All algorithms are written in Fortran 90/95 and compiled with GNU Fortran with flag -O3. 
The serial simulations are performed with a single CPU (Intel(R) Xeon(R) CPU E5440 @ 2.83GHz with 2G Memory) on an 8-core workstation. 
MPI parallel simulations are conducted 
on the DMC cluster (Intel(R) Xeon(R) CPU E5520 @ 2.27GHz with 1.5G memory for each core) at Alabama Supercomputing Center.  
Before seeing the numerical results, we define order and errors.  

\subsection{Order and Errors}
In this paper, we report the relative errors defined as
\begin{eqnarray}\label{eqerr1}
  {e_{\phi}}&=& \frac{\max\limits_{i=1,...,N}|\phi^{num}(x_{i})-\phi^{exa}(x_{i})|}{\max\limits_{i=1,...,N}|\phi^{exa}(x_{i})|}\label{eqerror}
\end{eqnarray}
where $N$ is the number of unknowns. 
Note $N$ is the number of vertices for HOBI-PB, 
the number of triangular elements for LOBI-PB,  
and the number of close-to-surface mesh points (irregular points) for APBS and MIBPB. 
The notation $\phi^{num}$ represents numerically solved surface potentials 
and $\phi^{exa}$ denotes the analytical solutions obtained 
by Kirkwood's spherical harmonic expansion \cite{Kirkwood34}. 
The discretization of APBS and MIBPB are on the Cartesian grid with mesh size $h$. 
The discretization of HOBI-PB and LOBI-PB are on the molecular surface with density $d$, 
number of vertices per \AA$^2$.

The numerical order of accuracy is computed with 
\begin{eqnarray}\label{eqeorder}
  \texttt{order}=\displaystyle\log_{\frac{\texttt{coarse\_mesh}}{\texttt{fine\_mesh}}}\frac{\texttt{coarse\_error}}{\texttt{fine\_error}}
\end{eqnarray}
following the convention of numerical analysis, where ``mesh" refers to $h$ for finite difference methods 
or density $d$ for boundary integral methods,  both at coarse and fine levels. 

\subsection{On a Spherical Cavity with one Centered Charge}
We first solve the linear Poisson-Boltzmann equation on a spherical cavity with radius 2\AA~and a centered change 1$e_c$ submerged in water with zero ionic strength. 
We report the electrostatic solvation energy $E_{sol}$ and surface potential errors $e_{\phi}$ computed with all above-mentioned methods in Table \ref{tb1}. The results of APBS and MIBPB are from reference \cite{Gengjcp07}.

\begin{table}
  \small{
    \caption{\small Accuracy comparison of different PB solvers on a spherical cavity (radius=2\AA~, $q(0,0,0)=1e_c$, $\veps=80$, $\kappa=0$); $h$ the mesh size; $d$ the number of vertices per~\AA$^2$.}\label{tb1}
  \begin{tabular}{c|ccc|ccc||c|ccc|ccc}
    \hline
     & \multicolumn{3}{c|}{APBS} & \multicolumn{3}{c||}{MIBPB} & & \multicolumn{3}{c|}{HOBI-PB} & \multicolumn{3}{c}{LOBI-PB}\\\hline
$h$  &$E_{sol}$ &$e_{\phi}$  &ord.&$E_{sol}$ &$e_{\phi}$&ord.   & $d$    &$E_{sol}$&$e_{\phi}$&ord. &$E_{sol}$ &$e_{\phi}$& ord.\\\hline
1	 &-83.44	&1.94e+0	&     &-81.95	&1.24e-2	&     &5	&-81.98	&1.45e-4 &            &-83.41&4.07e-4&   \\
0.5	 &-85.85	&1.31e+0	&0.6&-81.98	&1.91e-3	&2.7&10	&-81.98	&5.26e-5 & 1.5      &-83.13&2.55e-4&0.7\\
0.2	 &-82.58	&5.76e-1	&0.9&-81.98	&3.87e-4	&1.7&20	&-81.98	&1.52e-5 & 1.8      &-82.82&1.82e-4&0.5\\
0.1	 &-82.27	&2.94e-1	&1.0&-81.98	&1.07e-4	&1.9&40	&-81.98	&6.13e-6 & 1.3      &-82.60&1.27e-4&0.5\\
0.05 &-82.03	&1.49e-1	&1.0&-81.98	&2.31e-5	&2.2&80	&-81.98	&1.85e-6 & 1.7      &-82.43&8.58e-5&0.6\\
    \hline
  \end{tabular}}
\end{table}

From Table~{\ref{tb1}}, we can see that APBS provides acceptable value 
in electrostatic solvation energy compared with the true value 81.98 kcal/mol 
since only the potential at the center of the spherical cavity is required 
to compute the electrostatic solvation energy. 
However, the surface potentials computed by APBS method show large errors. 
This is due to the approximation on interface conditions and singular charges
of standard finite difference methods. 

MIBPB uses a more sophisticated  finite difference scheme.
The rigorous treatment on interface conditions 
and singular charges significantly improves accuracy. 
The electrostatic solvation energy is nearly perfect even at coarse grid 
and the surface potential is very accurate 
with solid 2nd order convergence pattern relative to mesh-size $h$.

LOBI-PB gives the electrostatic solvation energy in the same accurate level as APBS 
but produces much more accurate surface potentials than APBS does. 
These surface potentials converge at the 0.5th order relative to density $d$. 

HOBI-PB is obviously a more accurate method. 
It shows 1.5th order of convergence relative to density $d$  as reflected from surface potential errors 
and its accuracy is even better than results from MIBPB. 

These results demonstrate HOBI-PB is the most accurate PB solver among its peers. According to our knowledge, for this classic benchmark test, there is no other PB solvers can achieve 
the same level of accuracy as HOBI-PB does.

\subsection{A Spherical Cavity with One Eccentric Charge}

We further investigate the performance of HOBI-PB and LOBI-PB on a spherical cavity with radius 1\AA ~ and a 1$e_c$ eccentric charge at different locations. The errors shown are in terms of surface potential $e_u$ as previously defined. The Debye--H\"{u}ckel parameter $\kappa$ is set to 1.
The charge moves from the center of the sphere toward to the surface 
and we test the performance of the PB solvers 
in response to the change of locations. 
The closer the charge is to the surface, the more variant the induced charges on the surface appear to be. 
This interesting phenomenon draws attention from many researchers, e.g. the peak separation method by Juffer et al.~\cite{JBKPB} and the image method by Deng et al.~\cite{DengCaijcp2007}. In addition to that, we also try to investigate if the quality of the triangulation 
will affect the accuracy. 
To this end, we report the results computed on the triangular surfaces 
generated by MSMS \cite{Sanner}  and by 
a uniform icosahedral triangulation  routine (ico) \cite{Baumgardner}. 
MSMS generates triangles in various shapes and sizes, 
while ico generates uniform and equilateral triangles 
with fixed numbers of triangles such as 20, 80, 320, etc. 
The results are plotted in Fig.~\ref{fig_poterr}, 
and we have the following observation and discussion. 

\begin{figure}
\flushleft
(a)\hs{8.2cm}(b)\hs{1cm}\\
\vskip -2pt
\includegraphics[width=2.9in]{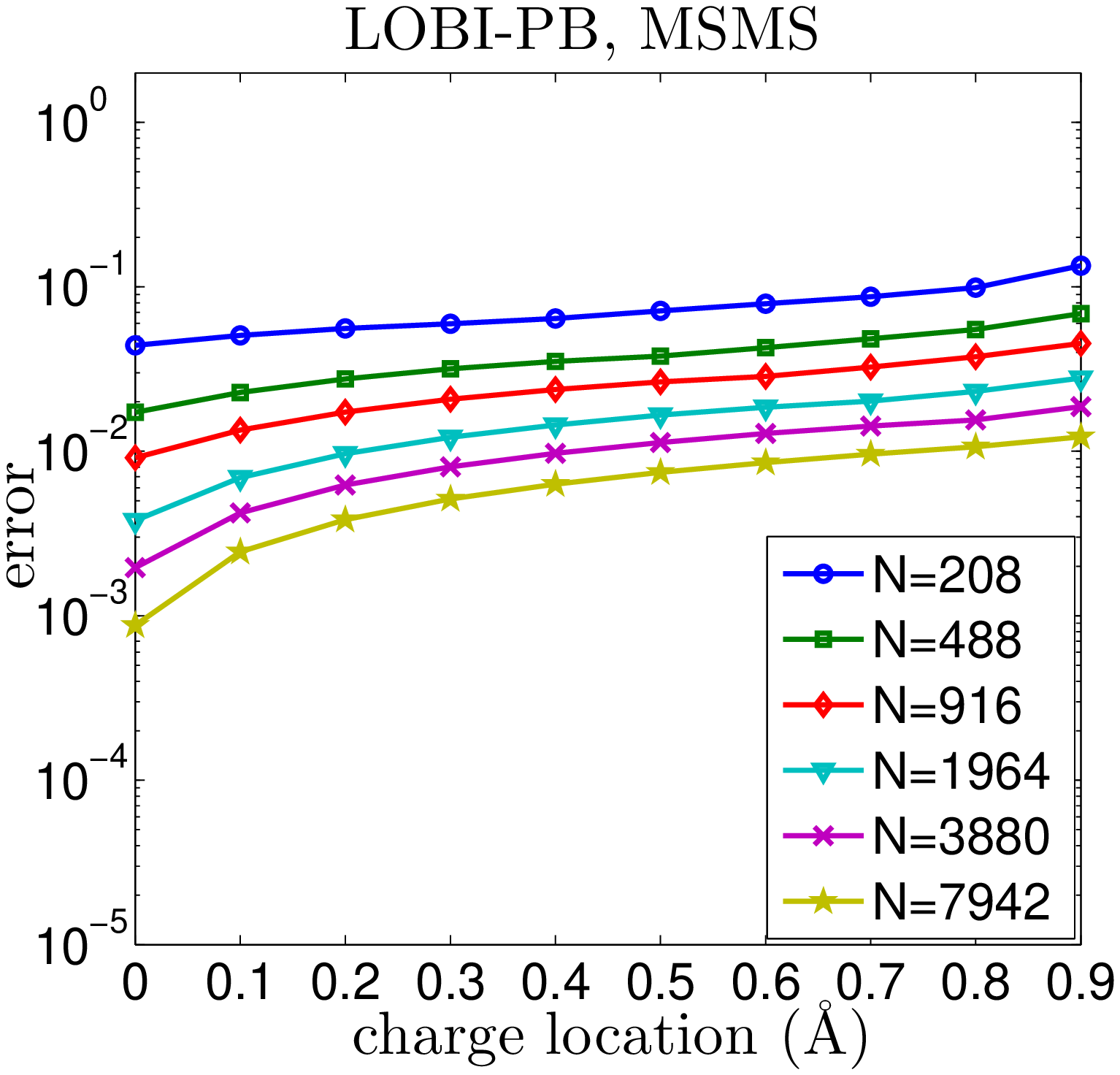}
\hs{1cm}
\includegraphics[width=2.9in]{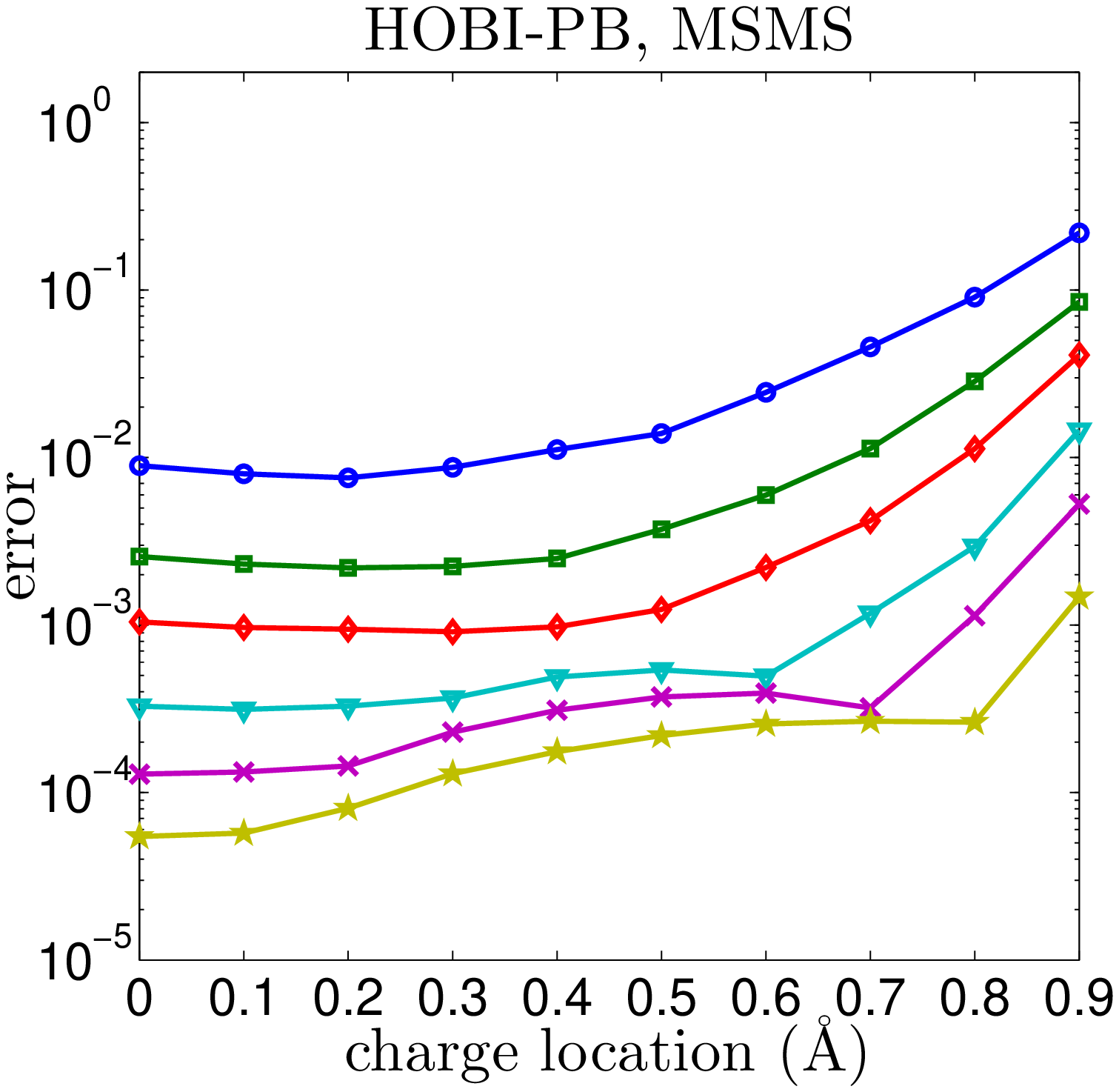}\\
(c)\hs{8.2cm}(d)\hs{1cm}\\
\vskip -2pt
\includegraphics[width=2.9in]{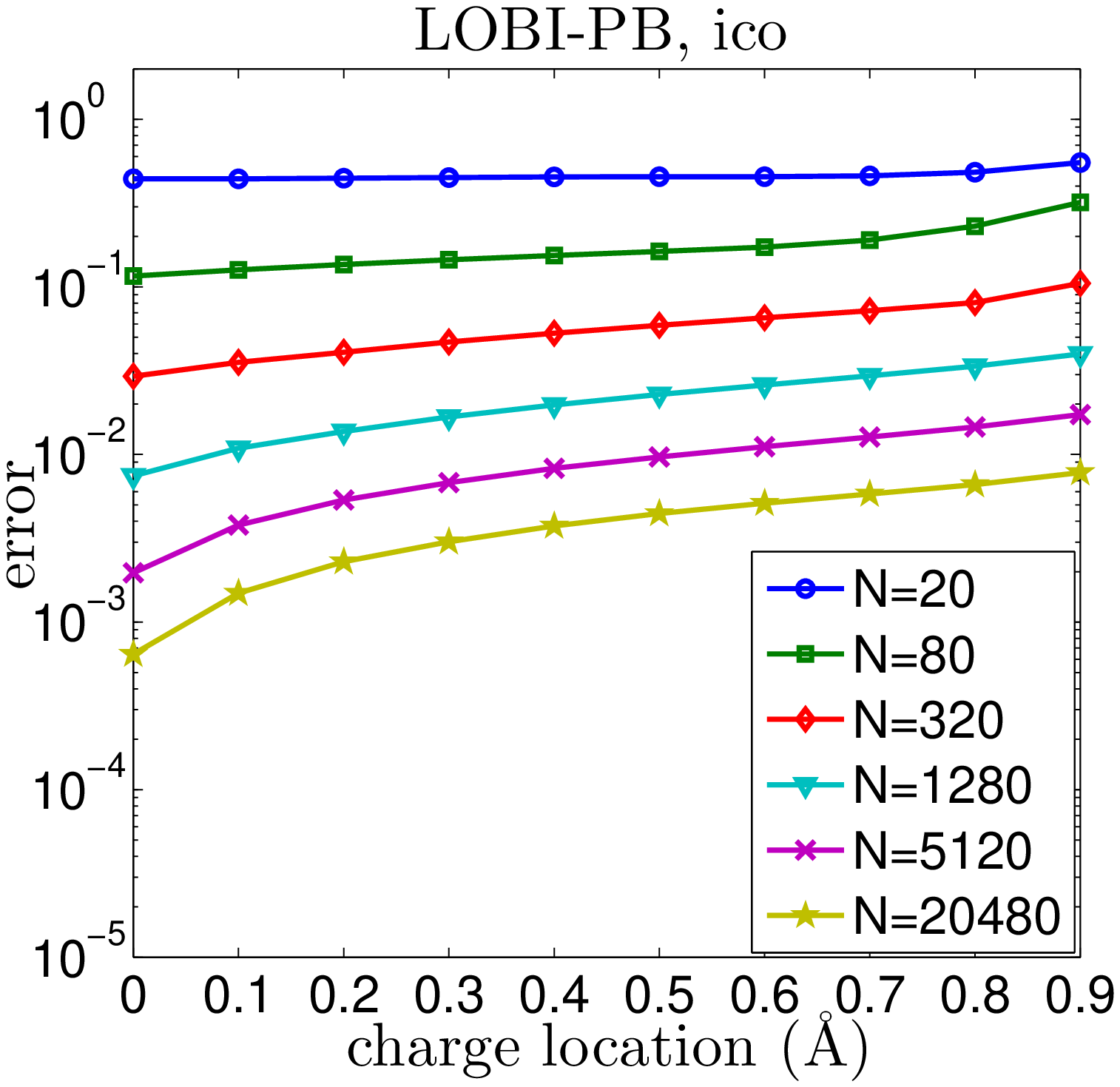}
\hs{1cm}
\includegraphics[width=2.9in]{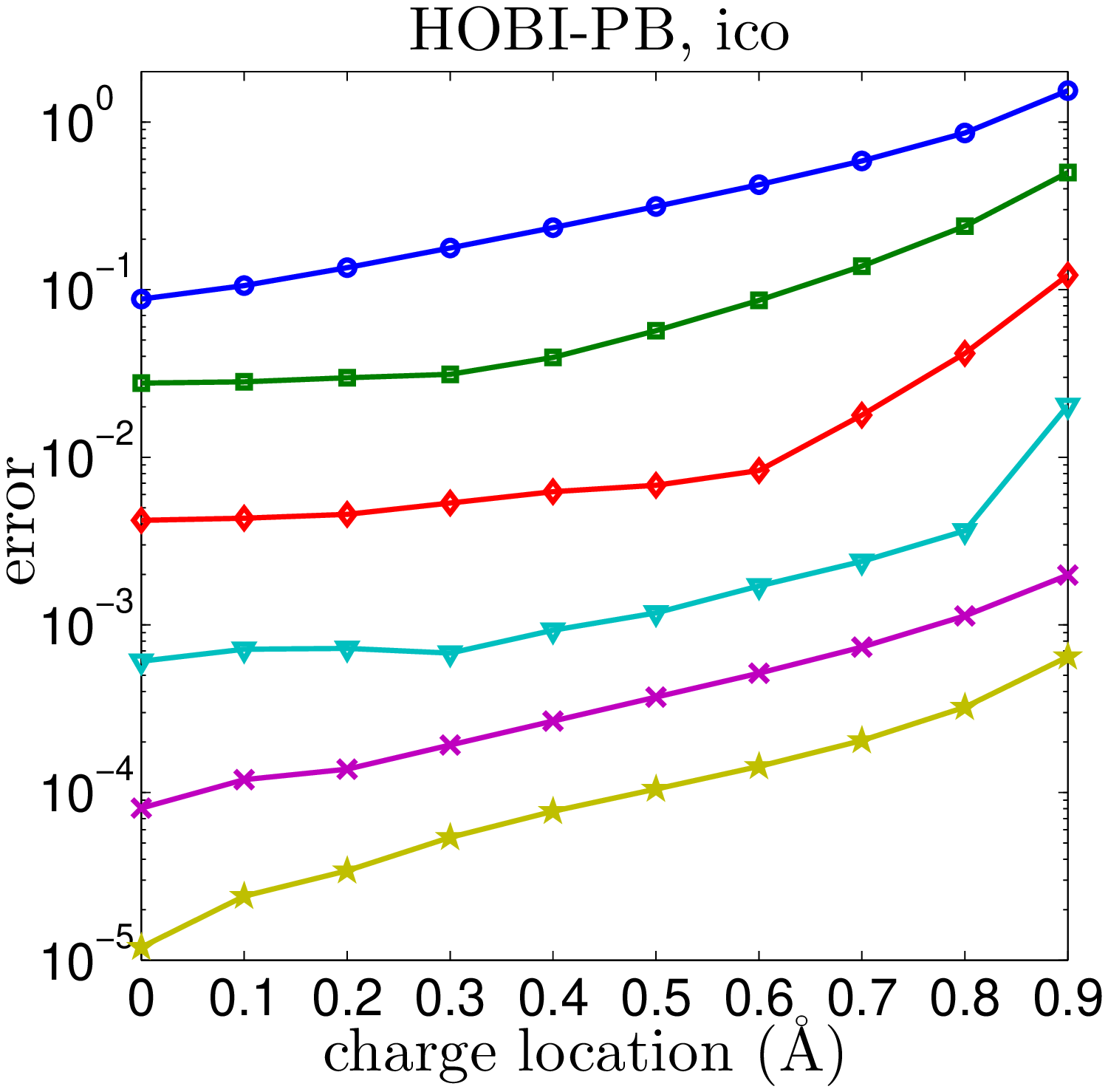}\\
(e)\hs{8.2cm}(f)\hs{1cm}\\
\vskip -2pt
\includegraphics[width=2.9in]{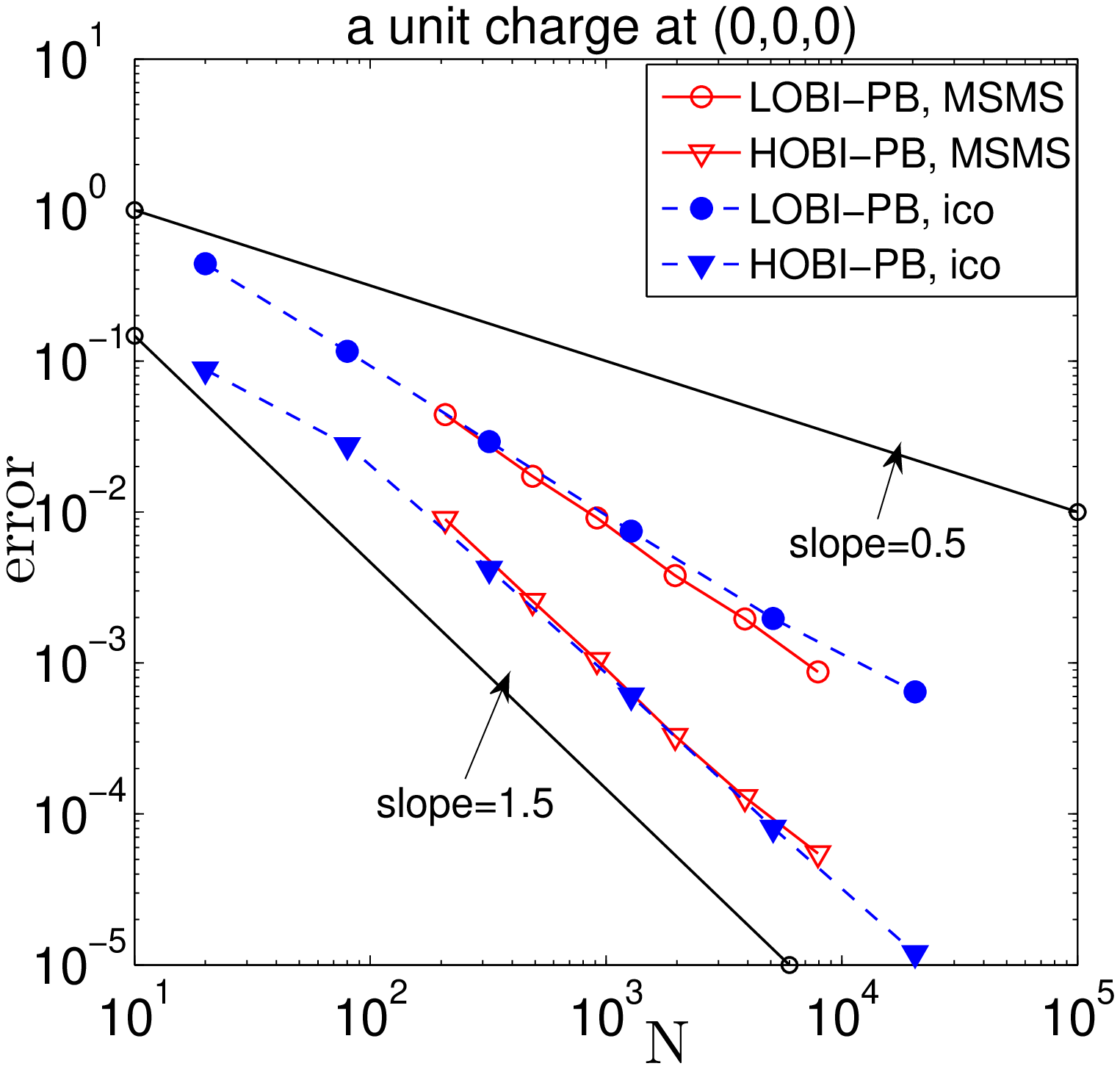}
\hs{1cm}
\includegraphics[width=2.9in]{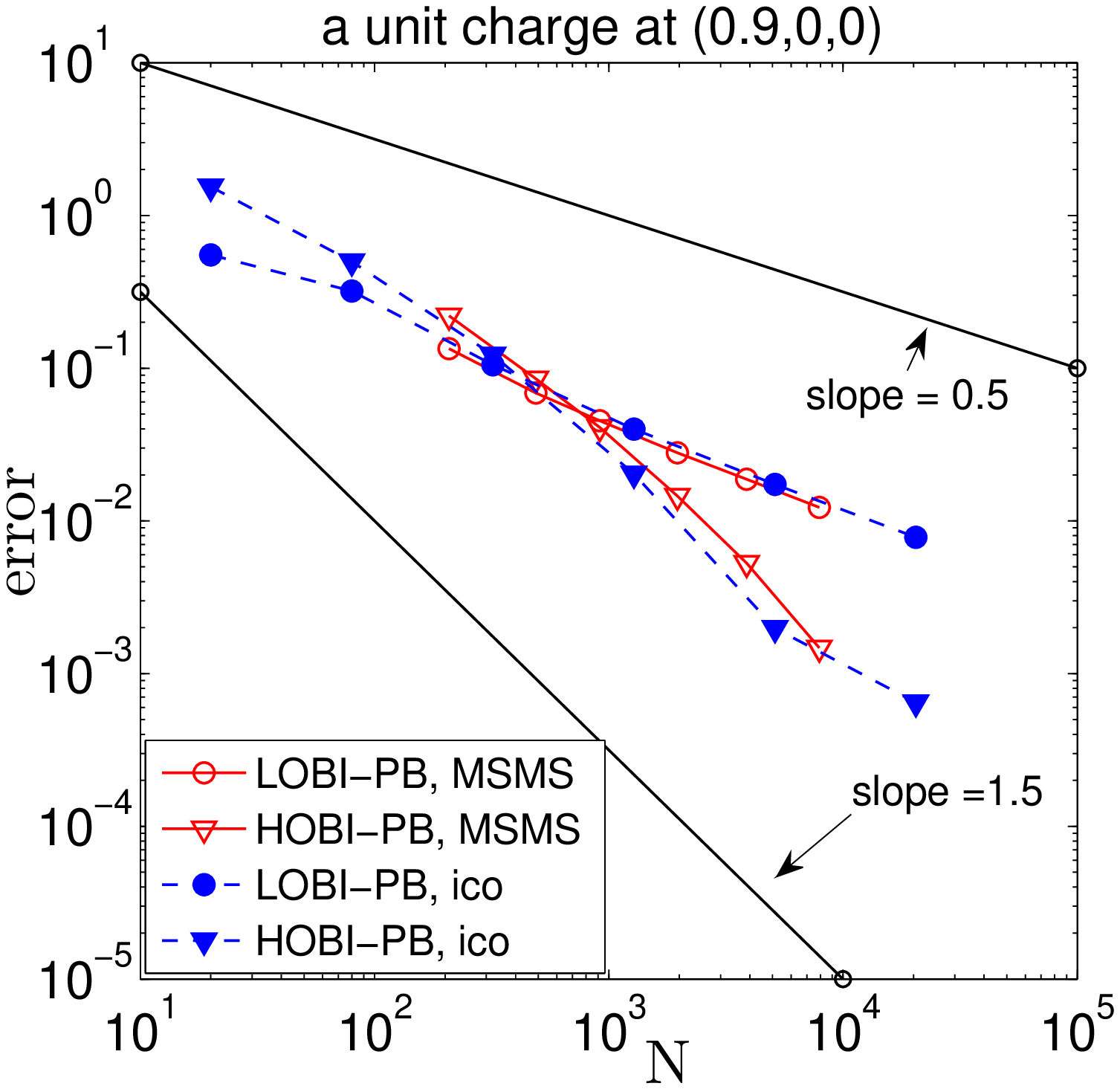}\\
\vskip -5pt
\caption{\small relative surface potential errors ($e_{\phi}$) on a spherical cavity with radius 1\AA~ and eccentric unit charge, $\kappa=1$, and $\varepsilon=80$: (a) LOBI-PB, MSMS; (b) HOBI-PB, MSMS;  (c) LOBI-PB, ico; (d) HOBI-PB, ico; (e) error vs. elements \# $N$, a charge located at (0,0,0); (f) error vs. elements \# $N$, a charge located at (0.9,0,0).}\label{fig_poterr}
\end{figure}

\noindent
(1) The results of HOBI-PB in Fig.~\label{fig_poterr}(b)(d) show significant improvement in accuracy compared with results of  LOBI-PB in Fig.~\label{fig_poterr}(a)(c). For HOBI-PB, the errors are smaller in general and the difference between two different meshes are larger, indicating better accuracy and higher-order convergence\\
(2) LOBI-PB, due to its simplicity in algorithm, 
shows more consistent convergence pattern. 
HOBI-PB, however affected by the quadrature rule and singularity removal schemes, 
shows some irregular patterns.\\ 
(3) We also observe that when the charge is located close to the interface e.g.~at (0.9,0,0),
for coarse mesh, LOBI-PB shows better accuracy than that from HOBI-PB. 
This is due to the fact that the close-to-boundary charge brings high variation of the induced charges on surface and multiple quadrature points selected in HOBI-PB amplify the variations. 
In LOBI-PB, there is only one quadrature point in each triangle at the centroid
thus the scheme is less sensitive to the induced charges.\\
(4) The quality of triangular surface will slightly affect the convergence of the boundary integral PB solvers. 
The results from icosahedral triangulation  routines show more uniform convergence pattern 
than that from MSMS.\\
(5) We further draw the error-vs-element plots for one charge located at (0,0,0) in Fig.~\ref{fig_poterr}(e) and (0.9,0,0) in Fig.~\ref{fig_poterr}(f) for both LOBI-PB (circle) and HOBI-PB (triangle) solvers on both MSMS (red, empty marker, solid line) and icosahedron (blue, filled marker, dashed line) triangles. By observing one color at a time (one kind of mesh at a time), we can see the pattern in Fig.~\ref{fig_poterr}(e) is uniform and HOBI-PB shows better accuracy (smaller $y$ values) and faster convergence (larger slope). The pattern in Fig.~\ref{fig_poterr}(f) is tangled. 
By observing the slope, we can still see that HOBI-PB converges faster than LOBI-PB generally. 
However, the error of the HOBI-PB is bigger than LOBI-PB for coarse mesh due to the interaction between the induced charges on surface and the singular charge near the surface, as explained in observation and discussion (3). In practice, the partial charges are at least 1-2~\AA~away from the molecular surfaces therefore the slightly deteriorated pattern observed here when charges are close to the surfaces will unlikely happen. \\

In short, the numerical results of HOBI-PB and some other reference PB solvers on the spherical cavities compared with available analytical solutions quantitatively demonstrate the achieved better accuracy and higher-order convergence of HIBO-PB. The complexity of the higher-order schemes for quadrature and regularization 
introduces only minor instabilities but gain significantly improved accuracy and convergence.   
Next we use HOBI-PB and LOBI-PB to solve PB equation and compute electrostatic solvation energy on a protein.  

\subsection{Computing Electrostatic Solvation Energy on Protein 1ajj}
The ultimate goal of HOBI-PB is to provide accurate electrostatic potentials for solvated biomolecules. 
With this solver, we solve the PB equation and compute the electrostatic solvation energy for many small-middle sized proteins. Here we take protein 1ajj with triangulation  at different resolutions as an example. The coordinates and partial charges of the protein are obtained from CHARMM22 force field \cite{MacKerell-CHARMM22}. 
The numerical results show that HOBI-PB solves the PB equation with significantly improved accuracy 
compared with LOBI-PB does. 

\begin{table}[htdp]
\caption{\small Testing results for protein 1ajj.
HOBI-PB results
showing 
electrostatic solvation energy $E_{sol}$, 
CPU time, 
memory usage;
number of GMRES iterations (it.).}
\begin{center}
\begin{tabular}{cccccc}
\hline
$d$	&\# of ele. $N$	&$E_{sol}$ (kcal/mol) 	&CPU (s)		&Memory (Mbyte)	&\# of it.\\\hline
1	&6027		&-1168.87		&115		&47	&10\\
2	&9198		&-1152.42		&266		&70	&10\\
4	&17278	&-1145.50		&867		&129	&9\\
8	&32386	&-1140.60		&3114		&238	&9\\
16	&66558	&-1139.23		&17790	&523	&13\\
32	&132028	&-1138.38		&56039	&759	&10\\
64	&270680	&-1138.49		&254198	&2116	&11\\\hline
\end{tabular}
\end{center}
\label{tb_1ajj}
\end{table}%

The numerical results for computing electrostatic solvation energy for protein 1ajj are reported in Table~\ref{tb_1ajj}. To produce these results, we discretize the molecular surface of protein 1ajj at different densities as seen in the first column of the table. 
Different densities result in different numbers of triangular elements 
which are listed in the second column of the table, 
characterizing the dimension of problem. 
In the third column, we report electrostatic solvation energies. 
We can see these values are very close to each other at different resolutions 
and they converge to a value near about $-1138.49$ kcal/mol. 
In column 4, we report the CPU time and it increases at the order of $\mathcal{O}(N^2)$.
This reveals currently the most critical limitation of HOBI-PB 
as the $\mathcal{O}(N^2)$ 
computational cost eventually will make the HOBI-PB prohibitively expensive. 
To alleviate the pain, we take advantage of the convenient parallelization of boundary integral formulation, 
and we will see the parallelization performance next. 
Memory uses are shown in column 5 and we see a $\mathcal{O}(N)$ pattern, 
which is advantageous compared with the 3D mesh-based methods whose memory uses
increase at  $\mathcal{O}(N^2)$ or even $\mathcal{O}(N^3)$, 
where $N$ is the number of unknowns.  
The last column is the number of iterations, and these stable results 
attribute to the well-posed integral formulation in Eqs.~(\ref{eqbim_3}) and~(\ref{eqbim_4}). 

\begin{figure}
\begin{center}
\includegraphics[width=4.5in]{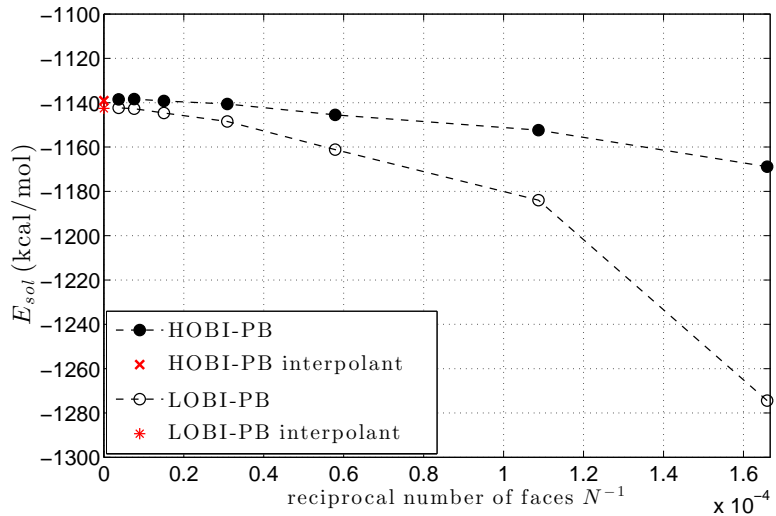}
\caption{\small Comparison of electrostatic solvation energy computed with HOBI-PB and LOBI-PB for protein 1ajj}\label{fig_1ajj}
\end{center}
\end{figure}

We next plot the solvation energies of HOBI-PB (solid circle) from Table~\ref{tb_1ajj} on Fig.~\ref{fig_1ajj} 
together with the solvation energies computed with LOBI-PB (empty circle) for the same protein at different resolutions.  The plot shows results from LOBI-PB converge to that of HOBI-PB. By using cubic interpolation, we can see both methods eventually converge toward almost the same values (the red ``x" for HOBI-PB and the red ``$*$" for LOBI-PB). The advantage of HOBI-PB is that it can achieve high accuracy even at very coarse mesh. 
For example, the first solid point from right, which corresponds to electrostatic solvation energy of $-1168.87$ kcal/mol at density $d=1$, is only $30$ kcal/mol different from the interpolated true value at about $-1139.09$ kcal/mol. Similar patterns are observed on our other tests on different proteins. 

We finally provided the MPI based parallel performance in terms of parallel efficiency of HOBI-PB for protein 1ajj at different meshes in Table~\ref{tb_1ajj_mpi}. Due to the limitation of the resources, our results are generated with up to 64 CPUs. The maximum memories for each CPU in MPI implementation is about 1.5G therefore for protein (1ajj) we solve the PB equation with maximum $d=32$ requiring 759M memory (note $d=64$ requires more than 2G memory per core from Table~\ref{tb_1ajj}).  From Table~\ref{tb_1ajj_mpi}, we can see the CPU time is substantially reduced when the code is run in parallel on high performance computers. For example, for $d=32$ with 132,028 elements, the serial work requires more than half of a day (50457s), the parallel work with 64 CPUs produce results in about 15 minutes  (860s). We see high parallel efficiency from the $\frac{T_1}{pT_p}$ column. When the dimension of the problem is sufficiently large (e.g. $N=66,558$ or $N=132,028$), the parallel efficiency with 64 CPUs is higher than 90\%. We observed occasionally the interesting larger-than-one parallel efficiencies and those could be explained by the traffic fluctuation on the cluster or possibly the argument mentioned in \cite{Parkinson}.  

\begin{table}[htdp]
\caption{\small MPI parallel performance for computing electrostatic solvation energy on protein 1ajj; $p$ is number of CPUs, $T_p$ is the time using $p$ CPUs, ${T_1}/{pT_p}$ is the parallel efficiency.}
\begin{center}
\begin{tabular}{ccccccc}
\hline
&\multicolumn{2}{c}{$N=132028$}& \multicolumn{2}{c}{$N=66558$} & \multicolumn{2}{c}{$N=32386$}\\\hline
$p$	& $T_p (s)$    & ${T_1}/{pT_p}$&	$T_p (s)$ & ${T_1}/{pT_p}$&$T_p (s)$ &${T_1}/{pT_p}$\\\hline	
1	&50457	&100.0\%	&17755	&100.0\%	&3060	&100.0\%\\
2	&24740	&102.0\%	&8790		&101.0\%	&1528	&100.1\%\\
4	&12877	&98.0\%	&4381		&101.3\%	&777	&98.5\%\\
8	&6398		&98.6\%	&2190		&101.3\%	&406	&94.3\%\\
16	&3321		&95.0\%	&1090		&101.8\%	&194	&98.8\%\\
32	&1677		&94.0\%	&571		&97.2\%		&106	&90.6\%\\
64	&860		&91.7\%	&306		&90.8\%		&62	&77.4\%\\\hline
\end{tabular}
\end{center}
\label{tb_1ajj_mpi}
\end{table}
In summary, the HOBI-PB solver solves the PB equation accurately on  both spherical cavities and real biomelecules. The surface potential and electrostatic solvation energy computed with the solver is accurate, fast-convergent and stable. The parallelization of the solver is easy to implement and the parallel efficiency is attractively high. 

\section{Conclusion}
This paper describes the schemes of a higher-order boundary integral
Poisson-Boltzmann (HOBI-PB) solver.  
This solver discretizes the molecular surface with curved triangles, and
performs numerical integral with four-point Gauss-Radau quadratures on regular triangles 
and with sixteen-point Gauss-Legendre quadratures on singular triangles. 
The singularities are regularized with a coordinate transformation. 
The numerical tests on spherical cavities show that HOBI-PB 
can achieve 1.5th order convergence of on surface potentials relative to boundary elements
and these computed surface potentials are very accurate even
at coarse mesh. In addition, the order of convergence
does not compromise when the electric charge is off-center or even closed to the
surface. The numerical tests of HOBI-PB on biomolecules show 
much improved accuracy compared with results from the popular LOBI-PB solvers.
The accurate surface potentials are of vital importance
to molecular modeling that are sensitive to electrostatics near or on the molecular surface.  
To improve the efficiency of HOBI-PB, we also developed its MPI based 
parallel version. The numerical results demonstrate very encouraging parallel efficiency, 
e.g. above 90\% when up to 64 CPUs work concurrently. 

HOBI-PB achieves higher accuracy at the price of more complex algorithms. 
The limitation of the HOBI-PB is mainly at the problem dimension it can treat, 
which is subject to the available computing resources. 
For example, for the accessibility  to up 64 CPUs each with 1.5G memory per core at Alabama Supercomputing Center,  the problem size HOBI-PB can handle for a reasonable long waiting time ($<$15 minutes) is about 150,000 elements. Considering the high accuracy, we can use fairly small density $1 \le d \le 5$, then PB equations on proteins with hundreds of residues can be conveniently solved.   
More computing resources will bring better performance on bigger problems. 
The rapid updating of computing power will definitely make HOBI-PB 
more and more capable. 

There are many spaces in which HOBI-PB can be improved and extended. For example, we are looking for better 
triangulation  programs for the molecular surfaces \cite{Bates08, Lujctc11, xu-zhang-09}. The currently adopted MSMS only provides 3 digits accuracy 
in vertices, positions and normal directions (we radially project vertices for spheres).
In addition, the adoption of fast algorithms such as FMM \cite{GH} and treecode \cite{Lijcp09} to HOBI-PB, although considerably 
challenging due to the complexity of HOBI-PB schemes, is under our consideration. A more challenging problem is the application of HOBI-PB to molecular dynamics \cite{Gengjcp11, Lu03}, where the PB equation will be solved at every time sampling. Furthermore, the higher-order schemes applied in HOBI-PB has the potential to solve other integral equations such as the integral forms of Helmholtz equation \cite{Amini} and Maxwell Equations \cite{Buffa}. For these challenges, the application of the quadrature rules and the treatment of singularity will be similar, however, new challenges such as obtaining the well-posedness of the integral formulation and applying the fast algorithms need further investigation.    
 
\section*{Acknowledgements}
The authors thank Robert Krasny, Yongcheng Zhou, Benzhuo Lu and Peijun Li for helpful discussions.
The work was supported by NSF grant DMS-0915057, University of Alabama new faculty startup fund 
and Alabama Supercomputer Center.


\end{document}